\documentclass[aap,MSNbibl,seceqn,citesort,dvips]{arximspdf}
\usepackage{graphicx}

\doi{10.1214/12-AAP859} 
\volume{23}
\issue{3}
\pubyear{2013}
\firstpage{895}
\lastpage{922}

\makeatletter

\newcommand{\cal}{\mathcal}

\newcommand{\EE}{\mathsf E}
\newcommand{\PP}{\mathsf P}
\newcommand{\cF}{{\cal F}}
\newcommand{\R}{\mathbb R}
\newcommand{\LL}{\mathbb L}

\newtheorem{lemma}{Lemma}
\newtheorem{theorem}[lemma]{Theorem}
\newtheorem{corollary}[lemma]{Corollary}
\newtheorem{proposition}[lemma]{Proposition}

\makeatother

\begin{document}
\begin{frontmatter}

\title{Three-dimensional Brownian motion and the~golden~ratio rule}
\runtitle{3D Brownian motion and the golden ratio rule}

\begin{aug}
\author[A]{\fnms{Kristoffer} \snm{Glover}\ead[label=e1]{kristoffer.glover@uts.edu.au}},
\author[A]{\fnms{Hardy} \snm{Hulley}\ead[label=e2]{hardy.hulley@uts.edu.au}}
\and
\author[B]{\fnms{Goran} \snm{Peskir}\corref{}\ead[label=e3]{goran@maths.man.ac.uk}}
\runauthor{K. Glover, H. Hulley and G. Peskir}
\affiliation{University of Technology, Sydney, University of Technology, Sydney and~The University of Manchester}
\address[A]{K. Glover \\
H. Hulley \\
Finance Discipline Group \\
University of Technology, Sydney \\
PO Box 123 \\
Broadway NSW 2007 \\
Australia \\
\printead{e1}\\
\hphantom{E-mail: }\printead*{e2}} 
\address[B]{G. Peskir \\
School of Mathematics \\
The University of Manchester \\
Oxford Road \\
Manchester M13 9PL \\
United Kingdom \\
\printead{e3}}
\end{aug}

\received{\smonth{1} \syear{2012}}
\revised{\smonth{3} \syear{2012}}

%
\begin{abstract}
Let $X=(X_t)_{t \ge0}$ be a transient diffusion process in
$(0,\infty)$ with the diffusion coefficient $\sigma>0$ and the scale
function $L$ such that $X_t \rightarrow\infty$ as $t \rightarrow
\infty$, let $I_t$ denote its running minimum for $t \ge0$, and let
$\theta$ denote the time of its ultimate minimum $I_\infty$. Setting
$c(i,x) = 1 - 2 L(x)/L(i)$ we show that the stopping time
\[
\tau_* = \inf\{ t \ge0 \vert X_t \ge f_*(I_t) \}
\]
minimizes $\EE(\vert\theta- \tau\vert- \theta)$ over all
stopping times $\tau$ of $X$ (with finite mean) where the optimal
boundary $f_*$ can be characterized as the minimal solution to
\[
f'(i) = - \frac{\sigma^2(f(i)) L'(f(i))}{c(i,f(i)) [L(f(i)) -
L(i)]} \int_i^{f(i)} \frac{c_i'(i,y) [L(y) - L(i)]}{\sigma^2(y)
L'(y)} \,dy
\]
staying strictly above the curve $h(i) = L^{-1}(L(i)/2)$ for $i>0$.
In particular, when $X$ is the radial part of three-dimensional
Brownian motion, we find that
\[
\tau_* = \inf\biggl\{ t \ge0 \Big\vert\frac{X_t - I_t}{I_t}
\ge\varphi\biggr\},
\]
where $\varphi= (1 + \sqrt{5})/2 = 1.61 \ldots$ is the golden ratio.
The derived results are applied to problems of optimal trading in the
presence of bubbles where we show that the golden ratio rule offers a
rigorous optimality argument for the choice of the well-known golden
retracement in technical analysis of asset prices.
\end{abstract}

%
\begin{keyword}[class=AMS]
\kwd[Primary ]{60G40}
\kwd{60J60}
\kwd{60J65}
\kwd[; secondary ]{34A34}
\kwd{49J40}
\kwd{60G44}.
\end{keyword}
\begin{keyword}
\kwd{Optimal prediction}
\kwd{transient diffusion}
\kwd{Bessel process}
\kwd{Brownian motion}
\kwd{the golden ratio}
\kwd{the maximality principle}
\kwd{Fibonacci retracement}
\kwd{support and resistance levels}
\kwd{constant elasticity of variance model}
\kwd{strict local martingale}
\kwd{bubbles}.
\end{keyword}

\end{frontmatter}

\section{Introduction}\label{sec1}

The \textit{golden ratio} has fascinated people of diverse interests
for at least 2400 years (see, e.g.,~\cite{Li}). In mathematics
(and the arts) two quantities $a$ and $b$ are in the golden ratio if
the ratio of the sum of the quantities $a + b$ to the larger quantity
$a$ is equal to the ratio of the larger quantity $a$ to the smaller
quantity~$b$. This amounts to setting\vspace*{1pt} $(a + b) / a = a/b
=: \varphi$ and solving $\varphi^2 - \varphi- 1 = 0$ which yields
$\varphi= (1 + \sqrt{5})/2 = 1.61 \ldots$ Apart from being abundant in
nature, and finding diverse applications ranging from architecture to
music, the golden ratio has also found more recent uses in
\textit{technical analysis} of asset prices (in strategies such as
\textit{Fibonacci retracement} representing an ad-hoc method for
determining \textit{support} and \textit{resistance} levels). Despite
its universal presence and canonical role in diverse applied areas, we
are not aware of any more \textit{exact} connections between the golden
ratio and \textit{stochastic processes} (including any proofs of
optimality in particular).\looseness=-1

One of the aims of the present paper is to disclose the appearance
of the golden ratio in an optimal stopping strategy related to the
radial part of three-dimensional Brownian motion. More specifically,
denoting the radial part by $X$ it is well known that $X$ is
transient in the sense that $X_t \rightarrow\infty$ as $t
\rightarrow\infty$. After starting at some $x>0$, the ultimate
minimum of $X$ will therefore be attained at some time $\theta$ that
is not predictable through the sequential observation of $X$ (in the
sense that it is only revealed at the end of time). The question we
are addressing is to determine a (predictable) stopping time of $X$
that is as close as possible to $\theta$. We answer this question by
showing that \textit{the time at which the excursion of $X$ away from
the running minimum $I$ and the running minimum $I$ itself form the
golden ratio} is as close as possible to $\theta$ in a normalized
mean deviation sense. We consider this problem by embedding it into
transient Bessel processes of dimension $d>2$ and in this context we
derive similar optimal stopping rules. We also disclose
further/deeper extensions of these results to transient diffusion
processes. The relevance of these questions in financial
applications is motivated by the problem of optimal trading in the
presence of bubbles. In this context we show that the golden ratio
rule offers a rigorous optimality argument for the choice of the
well-known \textit{golden retracement} in technical analysis of asset
prices. To our knowledge this is the first time that such an
argument has been found/given in the literature.

The problem considered in the present paper belongs to the class of
optimal prediction problems (within optimal stopping). Similar optimal
prediction problems have been studied in recent years by many authors
(see, e.g.,
\cite{BDP,Coh,DP-1,DP-2,DP-3,DPS,ET,GPS,NS,Ped-2,Sh-3,Sh-4,Ur}). Once
the ``unknown'' future is projected to the ``known'' present, we find
that the resulting optimal stopping problem takes a novel integral form
that has not been studied before. The appearance of the minimum process
in this context makes the problem related to optimal stopping problems
for the maximum process that were initially studied and solved in\vadjust{\goodbreak}
important special cases of diffusion processes in~\cite{DS,DSS} and
\cite{Ja}. The general solution to problems of this kind was derived in
the form of the maximality principle in~\cite{Pe-1}; see also
Section~13 and Chapter V in~\cite{PS} and the other references therein.
More recent contributions and studies of related problems include
\mbox{\cite{CHO,Ga-1,Ga-2,GZ,Ho,Ob-1,Ob-2,Ped-1}}. Close three-dimensional
relatives of these problems also appear in the recent papers~\cite{DGM}
and~\cite{Zi} where the problems were effectively solved by guessing
and finding the optimal stopping boundary in a closed form. The
maximality principle has been extended to three-dimensional problems in
the recent paper~\cite{Pe-4}.

Although the structure of the present problem is similar to some of
these problems, it turns out that none of these results is
applicable in the present setting. Governed by these particular
features in this paper we show how the problem can be solved when
(i) no closed-form solution for the candidate stopping boundary is
available and (ii) the loss function takes an integral form where
the integrand is a functional of both the process $X$ and its
running minimum $I$. This is done by extending the arguments
associated with the maximality principle to the setting of the
present problem and disclosing the general form of the solution that
is valid in all particular cases. The key novel ingredient revealed
in the solution is the replacement of the diagonal and its role in
the maximality principle by a nonlinear curve in the two-dimensional
state space of $X$ and $I$. We believe that this methodology is of
general interest and the arguments developed in the proof should be
applicable in similar two/multi-dimensional integral settings.

\section{Optimal prediction problem}\label{sec2}

1. We consider a nonnegative diffusion process $X=(X_t)_{t \ge0}$
solving
%
%
\begin{equation} \label{21}
dX_t = \mu(X_t) \,dt + \sigma(X_t) \,dB_t,
\end{equation}
where $\mu$ and $\sigma>0$ are continuous functions satisfying
(\ref{24}) and (\ref{25}) below, and $B=(B_t)_{t \ge0}$ is a standard
Brownian motion. By $\PP_{ x}$ we denote the probability measure
under which the process $X$ starts at $x > 0$. Recalling that the
scale function of $X$ is given by
%
%
\begin{equation} \label{22}
L(x) = \int^x \exp\biggl( - \int^y \frac{\mu(z)}{(\sigma^2
/2)(z)} \,dz \biggr) \,dy,
\end{equation}
and the speed measure of $X$ is given by
%
%
\begin{equation} \label{23}
m(dx) = \frac{dx}{(\sigma^2 /2)(x) L'(x)},
\end{equation}
we assume that the following conditions are satisfied:
%
%
\begin{eqnarray}
\label{24}
L(0+) &=& -\infty\quad\mbox{and}\quad L(\infty-) = 0, \\
\label{25}
\qquad\int_{0+}^1 L(dy) &=& \infty,\qquad
\int_{0+}^1 m(dy) < \infty\quad\mbox{and}\quad \int_{0+}^1
\vert L(y) \vert m(dy) < \infty.
\end{eqnarray}
From (\ref{24}) we read that $X$ is a transient diffusion process
in the sense that $X_t \rightarrow\infty$ $\PP_{ x}$-a.s. as $t
\rightarrow\infty$, and from (\ref{25})\vadjust{\goodbreak} we read that $0$ is an
entrance boundary point for $X$ in the sense that the process $X$
could start at $0$ but will never return to it (implying also that
$X$ will never visit $0$ after starting at $x>0$).\vspace*{8pt}

2. The main example we have in mind is the $d$-dimensional
Bessel process $X$ solving
%
%
\begin{equation} \label{26}
dX_t = \frac{d - 1}{2 X_t} \,dt + dB_t,
\end{equation}
where $d>2$. Recalling that the scale function is determined up to
an affine transformation we can choose the scale function
(\ref{22}) and hence the speed measure (\ref{23}) to read
%
%
\begin{eqnarray}
\label{27}
L(x) &=& - \frac{1}{x^{d-2}}, \\
\label{28}
m(dx) &=& \frac{2}
{d-2} x^{d-1} \,dx
\end{eqnarray}
for $x>0$. It is well known that when $d \in\{3,4, \ldots\}$ one
can realize $X$ as the radial part of $d$-dimensional standard
Brownian motion. Similar interpretations of (\ref{26}) are also
valid when $d=1$ (with an addition of the local time at zero) and
$d=2$ but $X$ is not transient in these cases (but recurrent), and
hence the problem considered below will have a trivial solution.
Other examples of (\ref{21}) are obtained by composing Bessel
processes solving (\ref{26}) with strictly decreasing and smooth
functions. This is of interest in financial applications and will be
discussed below. There are also many other examples of transient
diffusion processes solving (\ref{21}) that are not related to
Bessel processes.\vspace*{8pt}

3. To formulate the problem to be studied below consider the
diffusion process $X$ solving (\ref{21}), and introduce its running
minimum process $I=(I_t)_{t \ge0}$ by setting
%
%
\begin{equation} \label{29}
I_t = \inf_{0 \le s \le t} X_s
\end{equation}
for $t \ge0$. Due to the facts that $X$ is transient
(converging to $+\infty$) and $0$ is an entrance boundary point
for $X$, we see that the ultimate infimum $I_\infty= \inf_{ t \ge
0} X_t$ is attained at some random time $\theta$ in the sense that
%
%
\begin{equation} \label{210}
X_\theta= I_\infty
\end{equation}
with $\PP_{ x}$-probability one for $x>0$ given and fixed (the
case $x=0$ being trivial and therefore excluded). It is well
known that $\theta$ is unique up to a set of
$\PP_{ x}$-probability zero (cf.~\cite{Wi}, Theorem 2.4).
The random time $\theta$ is clearly unknown at any given time and
cannot be detected through sequential observations of the sample
path $t \mapsto X_t$ for $t \ge0$. In many applied situations of
this kind, we want to devise sequential strategies which will enable
us to come as ``close'' as possible to $\theta$. Most\vadjust{\goodbreak} notably, the
main example we have in mind is the problem of optimal trading in
the presence of bubbles to be addressed below. In mathematical terms
this amounts to finding a stopping time of $X$ that is as ``close'' as
possible to $\theta$. A first step toward this goal is provided by
the following lemma. We recall that stopping times of $X$ refer to
stopping times with respect to the natural filtration of $X$ that is
defined by ${\cal F}_t^X = \sigma(X_s \vert0 \le s \le t)$ for
$t \ge0$.
\begin{lemma}\label{lemma1}
We have
%
%
\begin{equation} \label{211}
\vert\theta- \tau\vert= \theta+ \int_0^\tau\bigl( 2 I(\theta
\le t) - 1 \bigr) \,dt
\end{equation}
for all stopping (random) times $\tau$ of $X$.
\end{lemma}
\begin{pf}
The identity is well known (see, e.g.,
\cite{PS}, page 450) and can be derived by noting that
%
%
\begin{eqnarray} \label{212}
\vert\theta- \tau\vert
&=&
(\theta- \tau)^+ + (\tau- \theta)^+
= \int_0^\theta I(\tau\le t) \,dt + \int_0^\tau I(\theta\le t)
\,dt \nonumber\\
&=& \int_0^\theta\bigl( 1 - I(\tau> t) \bigr) \,dt +
\int_0^\tau I(\theta\le t) \,dt \nonumber\\
&=& \theta- \int_0^\tau
I(\theta> t) \,dt + \int_0^\tau I(\theta\le t) \,dt \\
&=& \theta- \int_0^\tau\bigl(1 - I(\theta\le t) \bigr) \,dt +
\int_0^\tau I(\theta\le t) \,dt \nonumber\\
&=& \theta+ \int_0^\tau
\bigl( 2 I(\theta\le t) - 1 \bigr) \,dt\nonumber
\end{eqnarray}
for all stopping (random) times $\tau$ of $X$ as claimed.
\end{pf}

4. Taking $\EE_x$ on both sides in (\ref{211}) yields a nontrivial
measure of error (from $\tau$ to $\theta$) as long as $\EE_x
\theta< \infty$ for $x>0$ given and fixed. The latter condition,
however, may not always be fulfilled. For example, when $X$ is a
transient Bessel process of dimension $d>2$ it is known (see
\cite{Sh}, Lemma 1) that $\PP_{ x}(\theta> t) \sim t^{-(d/2-1)}$
as $t \rightarrow\infty$. Hence we see that $\EE_x \theta=
\int_0^\infty\PP_{ x}(\theta> t) \,dt < \infty$ if and only if
$d/2 - 1 > 1$ or equivalently $d > 4$. It is clear from
(\ref{211}), however, that the pointwise minimization of the
Euclidean distance on the left-hand side is equivalent to the
pointwise minimization of the integral on the right-hand side. To
preserve the generality we therefore ``normalize'' $\vert\theta-
\tau\vert$ on the left-hand side by subtracting $\theta$ from it.
After taking $\EE_x$ on both sides of the resulting identity, we
obtain
%
%
\begin{equation} \label{213}
\EE_x ( \vert\theta- \tau\vert- \theta) = \EE_x
\int_0^\tau\bigl( 2 I(\theta\le t) - 1 \bigr) \,dt
\end{equation}
for all stopping times $\tau$ of $X$ (for which the right-hand side
is well defined). The optimal prediction problem therefore becomes
%
%
\begin{equation} \label{214}
V(x) = \inf_\tau\EE_x ( \vert\theta- \tau\vert- \theta),
\end{equation}
where the infimum is taken over all stopping times $\tau$ of $X$
(with finite mean) and $x>0$ is given and fixed. Note that the
problem (\ref{214}) is equivalent to the problem of minimizing
$\EE_x \vert\theta- \tau\vert$ over all stopping times $\tau$ of
$X$ (with finite mean) whenever $\EE_x \theta< \infty$. To tackle
the problem (\ref{214}) we first focus on the right-hand side in
(\ref{213}) above.
\begin{lemma}\label{lemma2}
We have
%
%
\begin{equation} \label{215}
\EE_x \int_0^\tau\bigl( 2 I(\theta\le t) - 1 \bigr) \,dt =
\EE_x \int_0^\tau\biggl( 1 - 2 \frac{L(X_t)}{L(I_t)} \biggr)
\,dt
\end{equation}
for all stopping times $\tau$ of $X$ (with finite mean) and all $x >
0$.
\end{lemma}
\begin{pf}
Using a well-known argument (see, e.g.,
\cite{PS}, page 450), we find that
%
%
\begin{eqnarray} \label{216}
\EE_x \int_0^\tau\bigl( 2 I(\theta\le t) - 1 \bigr) \,dt
&=&
\EE_x \int_0^\infty\bigl( 2 I(\theta\le t) - 1 \bigr) I(t
< \tau) \,dt \nonumber\\
&=& \int_0^\infty\EE_x \bigl( \EE_x \bigl[\bigl( 2
I(\theta\le t) - 1 \bigr) I(t < \tau) \vert{\cal F}_t^X
\bigr] \bigr) \,dt \nonumber\\
&=& \int_0^\infty\EE_x \bigl( I(t < \tau)
\EE_x \bigl[\bigl( 2 I(\theta\le t) - 1 \bigr) \vert{\cal F}_t^X
\bigr] \bigr) \,dt \\
&=& \EE_x \int_0^\tau\bigl( 2 \PP_{ x}
(\theta\le t \vert{\cal F}_t^X) - 1 \bigr) \,dt \nonumber\\
&=&
\EE_x \int_0^\tau\bigl( 1 - 2 \PP_{ x} (\theta> t \vert
{\cal F}_t^X) \bigr) \,dt\nonumber
\end{eqnarray}
for any stopping time $\tau$ of $X$ (with finite mean) and any $x >
0$ given and fixed. Setting $I^t = \inf_{ s \ge t} X_s$ and
recalling that $I_t = \inf_{ 0 \le s \le t} X_s$, we find by the
Markov property that
%
%
\begin{eqnarray} \label{217}
\PP_{ x} (\theta> t \vert{\cal F}_t^X)
&=& \PP_{ x}(I^t <
I_t \vert{\cal F}_t^X) = \PP_{ x}(I^t < i \vert{\cal
F}_t^X) \vert_{i=I_t} \nonumber\\[-8pt]\\[-8pt]
&=& \PP_{ x}(I_\infty
\circ\theta_t < i \vert{\cal F}_t^X) \vert_{i=I_t} =
\PP_{X_t}(I_\infty< i) \vert_{i=I_t}\nonumber
\end{eqnarray}
for $t>0$. To compute the latter probability we recall that $M:=
L(X)$ is a continuous local martingale and note that $I_\infty<
i$ if and only if $L(I_\infty) < L(i)$ where $L(I_\infty) =
\inf_{ t \ge0} L(X_t) = \inf_{ t \ge0} M_t$. This shows
that the set $\{ I_\infty< i \}$ coincides with the set $\{
\inf_{ t \ge0} M_t < L(i) \}$ which in turn can be
expressed as $\{ \sigma< \infty\}$ where $\sigma= \inf\{
t \ge0 \vert M_t < L(i) \}$. Taking $x \ge i$ we see that the
continuous local martingale $M^\sigma= (M_{\sigma\wedge t})_{t \ge
0}$ is bounded above by $0$ and bounded below by $L(i)$ with\vadjust{\goodbreak}
$M_0^\sigma= L(x)$ under $\PP_{ x}$. It follows therefore that
$M^\sigma$ is a uniformly integrable martingale and hence by the
optional sampling theorem we find that
%
%
\begin{eqnarray} \label{218}
L(x) &=& \EE_x M_\sigma= \EE_x [L(i) I(\sigma< \infty)] +
\EE_x[ M_\infty I(\sigma= \infty)] \nonumber\\[-8pt]\\[-8pt]
&=& L(i) \PP_{ x}
(\sigma< \infty)\nonumber
\end{eqnarray}
upon using that $M_\infty:= \lim_{ t \rightarrow\infty} M_t =
0$ $\PP_{ x}$-a.s. on $\{ \sigma= \infty\}$. Combining
(\ref{218}) with the previous conclusions we obtain
%
%
\begin{equation} \label{219}
\PP_{ x}(I_\infty< i) = \PP_{ x}(\sigma< \infty) =
\frac{L(x)}{L(i)}
\end{equation}
for $i \le x$ in $(0,\infty)$. From (\ref{217}) and (\ref{219}) we
see that
%
%
\begin{equation} \label{220}
\PP_{ x} (\theta> t \vert{\cal F}_t^X) = \frac{L(X_t)}{L(I_t)}
\end{equation}
for all $x > 0$ and $t \ge0$ (for the underlying
three-dimensional law, see~\cite{CFS}, Theorem~A). Inserting this
expression back into (\ref{216}) we obtain (\ref{215}) and the
proof is complete.
\end{pf}

5. From (\ref{213}) and (\ref{215}) we see that the problem
(\ref{214}) is equivalent to
%
%
\begin{equation} \label{221}
V(x) = \inf_\tau\EE_x \int_0^\tau\biggl( 1 - 2 \frac{L(X_t)}
{L(I_t)} \biggr) \,dt,
\end{equation}
where the infimum is taken over all stopping times $\tau$ of $X$
(with finite mean) and $x>0$ is given and fixed. Passing from the
initial diffusion process $X$ to the scaled diffusion process $L(X)$
we see that there is no loss of generality in assuming that $\mu=0$
in (\ref{21}) or equivalently that $L(x)=x$ for $x>0$ (with
$L(X_t) \rightarrow0$ $\PP_{ x}$-a.s. as $t \rightarrow
\infty$). Note that the time of the ultimate minimum $\theta$ is the
same for both $X$ and $L(X)$ since $L$ is strictly increasing. Note
also that $\tau$ is a stopping time of $X$ if and only if $\tau$ is
a stopping time of $L(X)$. To keep the track of the general formulas
throughout we will continue with considering the general case (when
$\mu$ is not necessarily zero and $L$ is not necessarily the
identity function). This problem will be tackled in the next section
below.\vspace*{8pt}

6. For future reference we recall that the infinitesimal generator
of $X$ equals
%
%
\begin{equation} \label{222}
\LL_X = \mu(x) \,\frac{\partial}{\partial x} + \frac{\sigma^2(x)}{2}\,
\frac{\partial^2}{\partial x^2}
\end{equation}
for $x>0$. Throughout we denote $\tau_a = \inf\{ t \ge0
\vert X_t=a \}$ and set $\tau_{a,b} = \tau_a \wedge\tau_b$ for
$a<b$ in $(0,\infty)$. It is well known that
%
%
\begin{equation} \label{223}\hspace*{28pt}
\PP_{ x} ( X_{\tau_{a,b}}=a ) = \frac{L(b) - L(x)}{L(b)
- L(a)} \quad\mbox{and}\quad \PP_{ x} ( X_{\tau_{a,b}}=b ) =
\frac{L(x) - L(a)} {L(b) - L(a)}
\end{equation}
for $a \le x \le b$ in $(0,\infty)$. The Green function of $X$ is
given by
%
%
\begin{eqnarray} \label{224}
G_{a,b}(x,y) &=& \frac{(L(b) - L(y)) (L(x) - L(a))}{L(b)
- L(a)} \qquad\mbox{if } a \le x \le y \le b \nonumber\\[-8pt]\\[-8pt]
&=&
\frac{(L(b) - L(x)) (L(y) - L(a))}{L(b) - L(a)}
\qquad\mbox{if } a \le y \le x \le b.\nonumber
\end{eqnarray}
If $f\dvtx (0,\infty) \rightarrow\R$ is a measurable function, then it
is well known that
%
%
\begin{equation} \label{225}
\EE_x \int_0^{\tau_{a,b}} f(X_t) \,dt = \int_a^b f(y)
G_{a,b}(x,y) m(dy)
\end{equation}
for $a \le x \le b$ in $(0,\infty)$. This identity holds in the
sense that if one of the integrals exists, so does the other one, and
they are equal.

\section{Optimal stopping problem}\label{sec3}

It was shown in the previous section that the optimal prediction
problem (\ref{214}) is equivalent to the optimal stopping
problem~(\ref{221}). The purpose of this section is to present the
solution to the latter problem. Using the fact that the two problems
are equivalent this also leads to the solution of the former problem.

In the setting of (\ref{21})--(\ref{25}) consider the optimal
stopping problem (\ref{221}). This problem is two-dimensional and
the underlying Markov process equals $(I,X)$. Setting $I_t^i = i
\wedge\inf_{ 0 \le s \le t} X_s$ for $t \ge0$ enables $(I,X)$
to start at $(i,x)$ under $\PP_{ x}$ for $i \le x$ in $(0,\infty)$,
and we will denote the resulting probability measure on the
canonical space by $\PP_{ i,x}$. Thus under $\PP_{ i,x}$ the
canonical process $(I,X)$ starts at $(i,x)$. The problem
(\ref{221}) then extends as follows:
%
%
\begin{equation} \label{31}
V(i,x) = \inf_\tau\EE_{i,x} \int_0^\tau\biggl( 1 - 2
\frac{L(X_t)} {L(I_t)} \biggr) \,dt
\end{equation}
for $i \le x$ in $(0,\infty)$ where the infimum is taken over all
stopping times $\tau$ of $X$ (with finite mean). In addition to
$\sigma$ and $L$ from (\ref{21}) and (\ref{22}) above, let us set
%
%
\begin{equation} \label{32}
c(i,x) = 1 - 2 \frac{L(x)}{L(i)}
\end{equation}
for $i \le x$ in $(0,\infty)$. The main result of this section may
then be stated as follows.
\begin{theorem}\label{theor3}
The optimal stopping time in problem
(\ref{31}) is given by
%
%
\begin{equation} \label{33}
\tau_* = \inf\{ t \ge0 \vert X_t \ge f_*(I_t) \},
\end{equation}
where the optimal boundary $f_*$ can be characterized as the minimal
solution to
%
%
\begin{equation} \label{34}\quad
f'(i) = - \frac{\sigma^2(f(i)) L'(f(i))}{c(i,f(i)) [L(f(i)) -
L(i)]} \int_i^{f(i)} \frac{c_i'(i,y) [L(y) - L(i)]}{\sigma^2(y)
L'(y)} \,dy
\end{equation}
staying strictly above the curve $h(i) = L^{-1}(L(i)/2)$ for $i>0$
(in the sense that if the minimal solution does not exist, then there
is no optimal stopping time). The value function is given by
%
%
\begin{equation} \label{35}
V(i,x) = - \int_x^{f_*(i)} \frac{c(i,y) [L(y) - L(x)]}{(\sigma^2
/2)(y) L'(y)} \,dy
\end{equation}
for $i \le x \le f_*(i)$ and $V(i,x)=0$ for $x \ge f_*(i)$ with
$i>0$.
\end{theorem}
\begin{pf} 1. It is evident from the integrand in (\ref{31})
that the excursions of $X$ away from the running minimum $I$ play a
key role in the analysis of the problem. In particular, recalling
definition (\ref{32}), we see from (\ref{31}) that the process
$(I,X)$ can never be optimally stopped in the set $C_0:= \{ (i,x)
\in S \vert c(i,x) < 0 \}$ where we let $S = \{ (i,x) \in
(0,\infty) \times(0,\infty) \vert i \le x \}$ denote the
state space of the process $(I,X)$. Indeed, if $(i,x) \in C_0$ is
given and fixed, then the first exit time of $(I,X)$ from a
sufficiently small ball with the centre at $(i,x)$ (on which $c$ is
strictly negative) will produce a value strictly smaller than $0$
(the value corresponding to stopping at once). Defining
%
%
\begin{equation} \label{36}
h(i) = L^{-1} \bigl( \tfrac{1}{2} L(i) \bigr)
\end{equation}
for $i>0$ we see that $c(i,x)<0$ for $x<h(i)$ and $c(i,x)>0$ for
$x>h(i)$ whenever $i \le x$ in $(0,\infty)$ are given and fixed.
Note that the mapping $i \mapsto h(i)$ is increasing and continuous
as well as that $h(i)>i$ for $i>0$ with $h(0+)=0$ and $h(+\infty) =
+\infty$. This shows that $C_0 = \{ (i,x) \in S \vert i \le x
< h(i) \}$. Note in particular that $C_0$ contains the diagonal
$\{ (i,x) \in S \vert i=x \}$ in the state space.\vspace*{8pt}

2. Before we formalize further conclusions, let us recall that the
general theory of optimal stopping for Markov processes (see
\cite{PS}, Chapter 1) implies that the continuation set in the
problem (\ref{31}) equals $C = \{ (i,x) \in S \vert
V(i,x)<0 \}$, and the stopping set equals $D = \{ (i,x) \in S
\vert V(i,x)=0 \}$. It means that the first entry time of
$(I,X)$ into $D$ is optimal in problem (\ref{31}) whenever well
defined. It follows therefore that $C_0$ is contained in $C$, and the
central question becomes to determine the remainder of the set $C$.
Since $X_t \rightarrow\infty$ $\PP_{ x}$-a.s. as $t \rightarrow
\infty$ it follows that $L(X_t) \rightarrow0$ $\PP_{ x}$-a.s. as
$t \rightarrow\infty$ so that the integrand in (\ref{31}) becomes
strictly positive eventually, and this reduces the incentive to
continue (given also that the ``favorable'' set $C_0$ becomes more
and more distant). This indicates that there should exist a point
$f(i)$ at or above which the process $X$ should be optimally stopped
under $\PP_{ i,x}$ where $i \le x$ in $(0,\infty)$ are given and
fixed. This yields the following candidate:
%
%
\begin{equation} \label{37}
\tau_f = \inf\{ t \ge0 \vert X_t \ge f(I_t) \}
\end{equation}
for an optimal stopping time in (\ref{31}) where the function $i
\mapsto f(i)$ is to be determined.\vadjust{\goodbreak}

3. \textit{Free-boundary problem.} To compute the value function $V$
and determine the optimal function $f$, we are led to formulate the
free-boundary problem
%
%
\begin{eqnarray}
\label{38}
(\LL_X V)(i,x) &=& - c(i,x) \qquad\mbox{for } i < x < f(i),
\\
\label{39} V_i'(i,x) \vert_{x=i+} &=& 0\qquad
\mbox{(normal reflection)},
\\
\label{310} V(i,x) \vert_{x=f(i)-} &=& 0\qquad
\mbox{(instantaneous stopping)},
\\
\label{311} V_x'(i,x) \vert_{x=f(i)-} &=& 0\qquad
\mbox{(smooth fit)}
\end{eqnarray}
for $i>0$, where $\LL_X$ is the infinitesimal generator of $X$ given
in (\ref{222}) above. For the rationale and further details
regarding free-boundary problems of this kind, we refer to
\cite{PS}, Section 13, and the references therein (we note, in
addition, that the condition of normal reflection (\ref{39}) dates
back to~\cite{GSG}).\vspace*{8pt}

4. \textit{Nonlinear differential equation.} To solve the
free-boundary problem (\ref{38})--(\ref{311}) consider the stopping
time $\tau_f$ defined in (\ref{37}) and (formally) the resulting
function
%
%
\begin{equation} \label{312}
V_f(i,x) = \EE_{i,x} \int_0^{\tau_f} c(I_t,X_t) \,dt
\end{equation}
for $i \le x \le f(i)$ in $(0,\infty)$ given and fixed. Applying the
strong Markov property of $(I,X)$ at $\tau_{i,f(i)} = \inf\{ t
\ge0 \vert X_t \notin(i,f(i)) \}$ and using
(\ref{223})--(\ref{225}), we find that
%
%
\begin{eqnarray} \label{313}
V_f(i,x)= V_f(i,i) \frac{L(f(i)) - L(x)}{L(f(i)) - L(i)} +
\int_i^{f(i)} c(i,y) G_{i,f(i)}(x,y) m(dy).\hspace*{-35pt}
\end{eqnarray}
It follows from (\ref{313}) that
%
%
\begin{eqnarray} \label{314}
V_f(i,i) &=& \frac{L(f(i)) - L(i)}{L(f(i)) - L(x)} V_f(i,x)\nonumber\\[-8pt]\\[-8pt]
&&{} -
\frac{L(f(i)) - L(i)}{L(f(i)) - L(x)} \int_i^{f(i)} c(i,y)
G_{i,f(i)}(x,y) m(dy).\nonumber
\end{eqnarray}
Using (\ref{310}) and (\ref{311}) we find after dividing and
multiplying with $x - f(i)$ that
%
%
\begin{equation} \label{315}
\lim_{x \uparrow f(i)} \frac{V_f(i,x)}{L(f(i)) - L(x)} = - \frac{1}
{L'(f(i))} \,\frac{\partial V_f}{\partial x}(i,x) \bigg\vert_{x=f(i)-}
= 0.
\end{equation}
Moreover, it is easily seen by (\ref{224}) that
%
%
\begin{eqnarray} \label{316}
&&
\lim_{x \uparrow f(i)} \frac{L(f(i)) - L(i)}{L(f(i)) - L(x)}
\int_i^{f(i)} c(i,y) G_{i,f(i)}(x,y) m(dy)\nonumber\\[-8pt]\\[-8pt]
&&\qquad = \int_i^{f(i)}
c(i,y) [L(y) - L(i)] m(dy).\nonumber
\end{eqnarray}
Combining (\ref{314})--(\ref{316}) we see that
%
%
\begin{equation} \label{317}
V_f(i,i) = - \int_i^{f(i)} c(i,y) [L(y) - L(i)] m(dy).
\end{equation}
Inserting this back into (\ref{313}) and using
(\ref{224}) and (\ref{225}), we conclude that
%
%
\begin{equation} \label{318}
V_f(i,x) = - \int_x^{f(i)} c(i,y) [L(y) - L(x)] m(dy)
\end{equation}
for $i \le x \le f(i)$ in $(0,\infty)$. Finally, using (\ref{39})
we find that
%
%
\begin{eqnarray} \label{319}
f'(i) &=& - \frac{\sigma^2(f(i)) L'(f(i))}{2 c(i,f(i)) [L(f(i)) -
L(i)]}\nonumber\\[-8pt]\\[-8pt]
&&\hspace*{6.5pt}{}\times\int_i^{f(i)} c_i'(i,y) [L(y) - L(i)] m(dy)\nonumber
\end{eqnarray}
for $i>0$. Recalling that $C_0$ is contained in $C$, we see that
there is no restriction to assume that each candidate function $f$
solving (\ref{319}) satisfies $f(i) \ge h(i)$ for all $i>0$. In
addition we will also show below that all points $(i,h(i))$ belong
to $C$ for $i>0$ so that (at least in principle) there would be no
restriction to assume that each candidate function $f$ solving
(\ref{319}) also satisfies $f(i)>h(i)$ for all $i>0$. These
candidate functions will be referred to as admissible. We will also
see below, however, that solutions to (\ref{319}) ``starting'' at $h$
play a crucial role in finding/describing the solution.

Summarizing the preceding considerations, we can conclude that to
each candidate function $f$ solving (\ref{319}), there corresponds
the function (\ref{318}) solving the free-boundary problem
(\ref{38})--(\ref{311}) as is easily verified by direct
calculation. Note, however, that this function does not necessarily
admit the stochastic representation (\ref{312}) (even though it was
formally derived from this representation). The central question
then becomes how to select the optimal boundary $f$ among all
admissible candidates solving (\ref{319}). To answer this question
we will invoke the subharmonic characterization of the value
function (see~\cite{PS}, Chapter 1) for the three-dimensional Markov
process $(I,X,A)$ where $A_t = \int_0^t c(I_s,X_s) \,ds$ for $t \ge
0$. Fuller details of this argument will become clearer as we
progress below. It should be noted that among all admissible
candidate functions solving (\ref{319}) only the optimal boundary
will have the power of securing the stochastic representation
(\ref{312}) for the corresponding function (\ref{318}). This is a
subtle point showing the full power of the method (as well as
disclosing limitations of the optimal stopping problem itself).\vspace*{8pt}

5. \textit{The minimal solution.} Motivated by the previous question
we note from (\ref{318}) that $f \mapsto V_f$ is decreasing over
admissible solutions to (\ref{319}). This suggests to select the
candidate function among admissible solutions to (\ref{319})\vadjust{\goodbreak} that
is as far as possible from $h$. The subharmonic characterization of
the value function suggests to proceed in the opposite direction, and
this is the lead that we will follow in the sequel.

To address the existence and uniqueness of solutions to
(\ref{319}), denote the right-hand side of (\ref{319}) by
$\Phi(i,f(i))$. From the general theory of nonlinear differential
equations, we know that if the direction field $(i,f) \mapsto
\Phi(i,f)$ is (locally) continuous and (locally) Lipschitz in the
second variable, then the equation (\ref{319}) admits a (locally)
unique solution. For instance, this will be the case if, along a
(local) continuity of $(i,f) \mapsto\Phi(i,f)$, we also have a
(local) continuity of $(i,f) \mapsto\Phi_f'(i,f)$. In particular,
we see from the structure of $\Phi$ that equation (\ref{319})
admits a (locally) unique solution whenever $x \mapsto\sigma^2(x)$
is (locally) continuously differentiable. It is important to realize
that the preceding arguments apply only away from $h$ since each
point $(i,h(i))$ is a singularity point of equation (\ref{319})
in the sense that $f'(i+)=\infty$ when $f(i+)=h(i)$ due to
$c(i,h(i))=0$ for $i>0$. In this case it is also important to note
that the preceding arguments can be applied to the equivalent
equation for the inverse of $i \mapsto f(i)$ since this singularity
gets removed (the derivative of the inverse being zero).


%
%
\begin{figure}[b]

\includegraphics{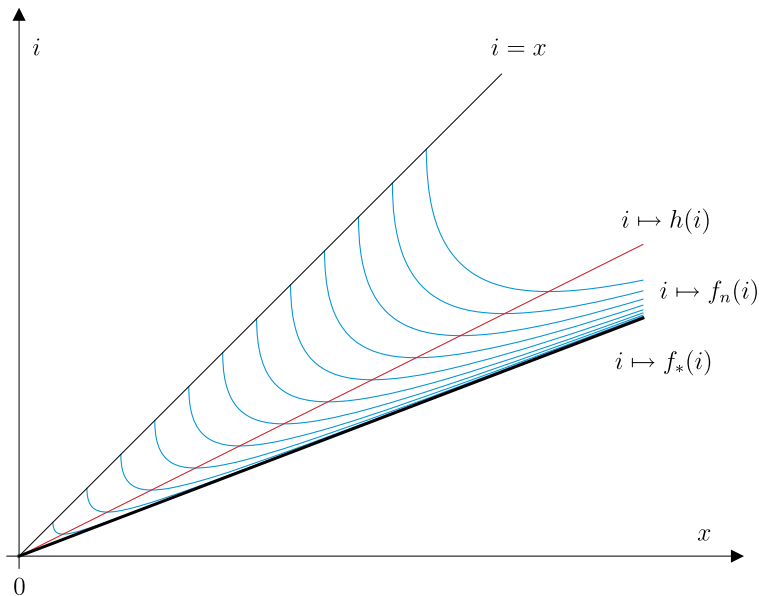}

\caption{Solutions $f_n$ and $f_*$ to the differential
equation (\protect\ref{34}) from Theorem \protect\ref{theor3}. The optimal stopping boundary
$f_*$ is the minimal solution staying strictly above the curve $h$.
This is a genuine drawing corresponding to the golden ratio rule of
Corollary \protect\ref{coro5} when $X$ is the radial part of three-dimensional
Brownian motion and the optimal stopping boundary $f_*$ is linear.}
\label{fig1}
\end{figure}


To construct the minimal solution to (\ref{319}) staying strictly
above $h$, we can proceed as follows (see Figure~\ref{fig1}).\vadjust{\goodbreak} For any
$i_n>0$ such that $i_n \downarrow0$ as $n \rightarrow\infty$, let
$i \mapsto f_n(i)$ denote the solution to (\ref{319}) on
$(i_n,\infty)$ such that $f_n(i_n+)=h(i_n)$. Note that $i \mapsto
f_n(i)$ is singular at $i_n$ and that passing to the equivalent
equation for the inverse of $i \mapsto f_n(i)$, this singularity gets
removed as explained above. (Note that the solution to the
equivalent equation for the inverse can be continued below $h(i_n)$
as well until hitting the diagonal at some strictly positive point
at which the derivative is $-\infty$. This yields another solution
to (\ref{319}) staying below $f_n$ and providing its ``physical''
link to the diagonal. We will not make use of this part of the
solution in the sequel.) Note that the right-hand side of
equation (\ref{319}) is positive for $f(i)>h(i)$ so that $i \mapsto
f_n(i)$ is strictly increasing on $[i_n,\infty)$. By the uniqueness
of the solution we know that the two curves $i \mapsto f_n(i)$ and
$i \mapsto f_m(i)$ cannot intersect for $n \ne m$, and hence we see
that $(f_n)_{n \ge1}$ is increasing. It follows therefore that $f_*
:= \lim_{ n \rightarrow\infty} f_n$ exists on $(0,\infty)$.
Passing to an integral equation equivalent to (\ref{319}) it is
easily verified that $i \mapsto f_*(i)$ solves (\ref{319}) wherever
finite. This $f_*$ represents the minimal solution to (\ref{319})
staying strictly above the curve $h$ on $(0,\infty)$. We will first
consider the case when $f_*$ is finite valued on $(0,\infty)$.\vspace*{8pt}

6. \textit{Stochastic representation.} We show that the function
(\ref{318}) associated with the minimal solution $f_*$ admits the
stochastic representation (\ref{312}). For this, let $i \mapsto
f_n(i)$ be the solution to (\ref{319}) on $(i_n,\infty)$ such that
$f_n(i_n+) = h(i_n)$ for $i_n>0$ with $i_n \downarrow0$ as $n
\rightarrow\infty$. Consider the function $(i,x) \mapsto
V_{f_n}(i,x)$ defined by (\ref{318}) for $i \le x \le f_n(i)$ and
$i \ge i_n$ with $n \ge1$ given and fixed. Recall that $V_{f_n}$
solves the free-boundary problem (\ref{38})--(\ref{311}) for $i \ge
i_n$. Consider the stopping time $\tau_n:= \tau_{i_n} \wedge
\tau_{f_n}$ where $\tau_{i_n} = \inf\{ t \ge0 \vert X_t=i_n \}$
and $\tau_{f_n} =\inf \{ t \ge0 \vert X_t \ge f_n(I_t) \}$.
Applying It\^o's formula and using (\ref{38}), we find that
%
%
\begin{eqnarray} \label{320}
V_{f_n}(I_{\tau_n},X_{\tau_n}) &=& V_{f_n}(i,x) + \int_0^{\tau_n}
\frac{\partial V_{f_n}}{\partial i}(I_t,X_t) \,dI_t\nonumber\\
&&{} + \int_0^{
\tau_n} \frac{\partial V_{f_n}}{\partial x}(I_t,X_t) \,dX_t\nonumber\\
&&{}+ \frac{1}{2} \int_0^{\tau_n} \frac{
\partial^2 V_{f_n}} {\partial x^2} (I_t,X_t) \,d\langle X,X\rangle_t
\nonumber\\[-8pt]\\[-8pt]
&=& V_{f_n}(i,x) + \int_0^{\tau_n} \LL_X(V_{f_n})
(I_t,X_t) \,dt\nonumber\\
&&{} + \int_0^{\tau_n} \sigma(X_t) \,\frac{\partial
V_{f_n}}{\partial x}(I_t,X_t) \,dB_t \nonumber\\
&=&
V_{f_n}(i,x) - \int_0^{ \tau_n} c(I_t,X_t) \,dt +M_{\tau_n},\nonumber
\end{eqnarray}
where we also use (\ref{39}) to conclude that the integral with
respect to $dI_t$ is equal to zero and $M_t = \int_0^{t \wedge
\tau_n} \sigma(X_s) (\partial V_{f_n}/\partial x)(I_s,X_s) \,dB_s$
is a continuous local martingale for $t \ge0$.

Since the process $(I,X)$ remains in the compact set $\{ (j,y) \in
S \vert i_n \le j \le y \le f_n(i) \}$ up to time $\tau_n$
under $\PP_{i,x}$, and both $\sigma$ and $\partial V_{f_n}/\partial
x$ are continuous (and thus bounded) on this set, we see that $M$ is
a uniformly integrable martingale, and hence by the optional sampling
theorem we have $\EE_{i,x}M_{\tau_n} = 0$. Taking $\EE_{i,x}$ on
both sides of (\ref{320}), we therefore obtain
%
%
\begin{eqnarray} \label{321}
V_{f_n}(i,x) &=& \EE_{i,x} V_{f_n}(I_{\tau_n},X_{\tau_n}) + \EE_{i,x}
\int_0^{ \tau_n} c(I_t,X_t) \,dt \nonumber\\[-8pt]\\[-8pt]
&=& V_{f_n}(i_n,i_n)\:
\PP_{i,x} (\tau_{i_n} < \tau_{f_n}) + \EE_{i,x} \int_0^{ \tau_n}
c(I_t,X_t) \,dt\nonumber
\end{eqnarray}
since $(I_{\tau_n},X_{\tau_n}) = (i_n,i_n)$ on $\{ \tau_{i_n} <
\tau_{f_n} \}$ and $V_{f_n}(I_{\tau_n},X_{\tau_n}) = 0$ on $\{
\tau_{f_n} < \tau_{i_n} \}$. Using that $\vert c \vert\le1$ we
find by (\ref{317}), (\ref{223}) and (\ref{36}) that
%
%
\begin{eqnarray} \label{322}
&&
\vert V_{f_n}(i_n,i_n) \vert\;\PP_{i,x} (\tau_{i_n} < \tau_{f_n})\nonumber\\
&&\qquad\le \int_{i_n}^{f(i_n)} \vert c(i,y) \vert\vert L(y) - L(i_n)
\vert m(dy) \;\PP_{i,x} (\tau_{i_n} < \tau_{f_n}) \nonumber\\[-8pt]\\[-8pt]
&&\qquad\le \vert L(h(i_n)) - L(i_n) \vert\int_{i_n}^{h(i_n)} m(dy)\:
\frac{L(f_*(i)) - L(x)}{L(f_*(i)) - L(i_n)} \nonumber\\
&&\qquad=
\frac{1}{2} \vert L(i_n) \vert\frac{L(f_*(i)) - L(x)}
{L(f_*(i)) - L(i_n)} \int_{i_n}^{h(i_n)} m(dy)
\rightarrow0\nonumber
\end{eqnarray}
as $n \rightarrow\infty$ since $L(i_n) \rightarrow-\infty$ and
$h(i_n) \rightarrow0$, so that $\int_{i_n}^{h(i_n)} m(dy)
\rightarrow0$, due to (\ref{25}) above. Hence letting $n
\rightarrow\infty$ in (\ref{321}) and using that $V_{f_n}
\rightarrow V_{f_*}$ by the monotone convergence theorem, as well as
that $\tau_n \uparrow\tau_{f_*}$ since $f_n \uparrow f_*$ and $i_n
\downarrow0$, we find noting that $\EE_{i,x} \tau_{f_*} < \infty$
and using the dominated convergence theorem that
%
%
\begin{equation} \label{323}
V_{f_*}(i,x) = \EE_{i,x} \int_0^{\tau_{f_*}} c(I_t,X_t) \,dt
\end{equation}
for all $i \le x $ in $(0,\infty)$ as claimed.\vspace*{8pt}

7. \textit{Nonpositivity.} We show that for every solution $f$ to
(\ref{319}) such that $f \ge f_*$ on $(0,\infty)$ and the function
$V_f$ defined by (\ref{318}) above, we have
%
%
\begin{equation} \label{324}
V_f(i,x) \le0
\end{equation}
for all $i \le x $ in $(0,\infty)$. Clearly, since $c(i,y) \ge0$
for $y \ge h(i)$ in (\ref{318}), it is enough to prove (\ref{324})
for $f_*$ and $i \le x < h(i)$ with $i>0$. For this, consider the
stopping time $\tau_h = \inf\{ t \ge0 \vert X_t \ge h(I_t) \}$
and note that $\tau_{f_*} = \tau_h + {\tau_{f_*}} \circ
{\theta_{\tau_h}}$. Hence by the strong Markov property of $(I,X)$
applied at $\tau_h$ we find using (\ref{323}) that
%
%
\begin{eqnarray} \label{325}
V_{f_*}(i,x) &=& \EE_{i,x} \int_0^{\tau_h} c(I_t,X_t) \,dt +
\EE_{i,x} \int_{\tau_h}^{\tau_h + \tau_{f_*} \circ\theta_{
\tau_h}} c(I_t,X_t) \,dt \nonumber\\
&=& \EE_{i,x} \int_0^{\tau_h}
c(I_t,X_t) \,dt + \EE_{i,x} \int_0^{\tau_{f_*} \circ\theta_{
\tau_h}} c(I_{t+\tau_h},X_{t+\tau_h}) \,dt \nonumber\\
&=& \EE_{i,x}
\int_0^{\tau_h} c(I_t,X_t) \,dt + \EE_{i,x} \EE_{i,x}\biggl[
\int_0^{\tau_{f_*}} c(I_t,X_t) \,dt \circ\theta_{\tau_h} \Big\vert
\cF_{\tau_h}^X \biggr] \\
&=& \EE_{i,x} \int_0^{\tau_h}
c(I_t,X_t) \,dt + \EE_{i,x} \EE_{I_{\tau_h},X_{\tau_h}}\biggl[
\int_0^{\tau_{f_*}} c(I_t,X_t) \,dt \biggr] \nonumber\\
&=& \EE_{i,x}
\int_0^{\tau_h} c(I_t,X_t) \,dt + \EE_{i,x} V_{f_*}(I_{\tau_h},
X_{\tau_h}) \le0,\nonumber
\end{eqnarray}
where the final inequality follows from the facts that $c(I_t,X_t)
\le0$ for all $t \in[0,\tau_h]$ and $V_{f_*}(I_{\tau_h},
X_{\tau_h}) \le0$ due to $X_{\tau_h} = h(I_{\tau_h})$ upon
recalling (\ref{318}) as already indicated above. This completes
the proof of (\ref{324}).\vspace*{8pt}

8. \textit{Optimality of the minimal solution.} We will begin by
disclosing the subharmonic characterization of the value function
(\ref{31}) in terms of the solutions to (\ref{319}) staying
strictly above $h$. For this, let $i \mapsto f(i)$ be any solution
to (\ref{319}) satisfying $f(i)>h(i)$ for all $i>0$. Consider the
function $(i,x) \mapsto V_f(i,x)$ defined by (\ref{318}) for $i \le
x \le f(i)$ in $(0,\infty)$ and set $V_f(i,x)=0$ for $x \ge f(i)$ in
$(0,\infty)$. Let $i \le x$ in $(0,\infty)$ be given and fixed. Due
to the ``double-deck'' structure of $V_f$, we can apply the
change-of-variable formula from~\cite{Pe-2} that in view of
(\ref{311}) reduces to standard It\^o's formula and gives
%
%
\begin{eqnarray} \label{326}
V_f(I_t,X_t) &=& V_f(i,x) + \int_0^t \frac{\partial V_f}{\partial i}
(I_s,X_s) \,dI_s \nonumber\\
&&{}+ \int_0^t \frac{\partial V_f}{\partial x}(I_s,X_s)
\,dX_s
+ \frac{1}{2} \int_0^t \frac{\partial^2 V_f}
{\partial x^2} (I_s,X_s) \,d\langle X,X\rangle_s \nonumber\\[-8pt]\\[-8pt]
&=& V_f(i,x) +
\int_0^t \LL_X(V_f)(I_s,X_s) \,ds \nonumber\\
&&{}+ \int_0^t \sigma(X_s)\,
\frac{\partial V_f}{\partial x}(I_s,X_s) \,dB_s,\nonumber
\end{eqnarray}
where we also use (\ref{39}) to conclude that the integral with
respect to $dI_s$ is equal to zero. The process $M=(M_t)_{t \ge0}$
defined by
%
%
\begin{equation} \label{327}
M_t = \int_0^t \sigma(X_s) \,\frac{\partial V_f}{\partial x}(I_s,X_s)
\,dB_s
\end{equation}
is a continuous local martingale. Introducing the increasing process
$P=(P_t)_{t \ge0}$ by setting
%
%
\begin{equation} \label{328}
P_t = \int_0^t c(I_s,X_s) I\bigl(X_s \ge f(I_s)\bigr) \,ds
\end{equation}
and using the fact that the set of all $s$ for which $X_s$ equals
$f(I_s)$ is of Lebesque measure zero, we see by (\ref{38}) that
(\ref{326}) can be rewritten as follows:
%
%
\begin{equation} \label{329}
V_f(I_t,X_t) + \int_0^t c(I_s,X_s) \,ds = V_f(i,x) + M_t + P_t.
\end{equation}
From this representation we see that the process $V_f(I_t,X_t) +
\int_0^t c(I_s,X_s) \,ds$ is a local submartingale for $t \ge0$.

Let $\tau$ be any stopping time of $X$ (with finite mean). Choose a
localization sequence $(\sigma_n)_{n \ge1}$ of bounded stopping
times for $M$. Then by (\ref{324}) and (\ref{329}) we can conclude
using the optional sampling theorem that
%
%
\begin{eqnarray} \label{330}\qquad
\EE_{i,x} \int_0^{\tau\wedge\sigma_n} c(I_t,X_t) \,dt
&\ge& \EE_{i,x} \biggl[ V_f(I_{\tau\wedge\sigma_n},X_{\tau\wedge
\sigma_n}) + \int_0^{\tau\wedge\sigma_n} c(I_t,X_t) \,dt
\biggr] \nonumber\\[-8pt]\\[-8pt]
&\ge&
V_f(i,x) + \EE_{i,x} M_{\tau\wedge
\sigma_n} = V_f(i,x).\nonumber
\end{eqnarray}
Letting $n \rightarrow\infty$ and using the dominated convergence
theorem (upon recalling that $\vert c \vert\le1$ as already used
above) we find that
%
%
\begin{equation} \label{331}
\EE_{i,x} \int_0^\tau c(I_t,X_t) \,dt \ge V_f(i,x).
\end{equation}
Taking first the infimum over all $\tau$, and then the supremum over
all $f$, we conclude that
%
%
\begin{equation} \label{332}
V(i,x) \ge\sup_f V_f(i,x) = V_{f_*}(i,x),
\end{equation}
upon recalling that $f \mapsto V_f$ is decreasing over $f \ge f_*$
so that the supremum is attained at $f_*$. Combining (\ref{332})
with (\ref{323}) we see that (\ref{33}) and (\ref{35}) hold as
claimed.

Note that (\ref{330}) implies that the function $(i,x) \mapsto
V_f(i,x) + a$ is subharmonic for the Markov process $(I,X,A)$ where
$A_t = \int_0^t c(I_s,X_s) \,ds$ for $t \ge0$. Recalling that $f
\mapsto V_f$ is decreasing over $f \ge f_*$, and that $V_f(i,x) \le
0$ for all $i \le x$ in $(0,\infty)$ by (\ref{324}) above, we see
that selecting the minimal solution $f_*$ staying strictly above $h$
is equivalent to invoking the subharmonic characterization of the
value function (according to which the value function is the largest
subharmonic function lying below the loss function). For more
details on the latter characterization in a general setting we refer
to~\cite{PS}, Chapter 1. It is also useful to know that the
subharmonic characterization of the value function represents the
dual problem\vadjust{\goodbreak} to the primal problem (\ref{31}) (for more details on
the meaning of this claim including connections to the Legendre
transform see~\cite{Pe-3}).

Consider finally the case when $f_*$ is not finite valued on
$(0,\infty)$. Since $i \mapsto f_*(i)$ is increasing we see that
there is $i_* \ge0$ such that $f_*(i)<\infty$ for all $i \in
(0,i_*)$ when $i_*>0$ and $f_*(i)=\infty$ for all $i \ge i_*$ with
$i \ne0$ when $i_*=0$. If $i_*>0$, then the proof above can be
applied in exactly the same way to show that (\ref{33}) and
(\ref{35}) hold as claimed under $\PP_{i,x}$ for all $i \le x$ in
$(0,\infty)$ with $i < i_*$. If $i \ge i_*$ with $i \ne0$ when
$i_*=0$, then the same proof shows that (\ref{35}) still holds with
$\infty$ in place of $f_*(i)$, however, the stopping time
(\ref{33}) can no longer be optimal in (\ref{31}). This is easily
seen by noting that the value in (\ref{35}) is nonpositive (it
could also be $-\infty$) for any $x \ge h(i)$ for instance, while
the $\PP_{i,x}$-probability for $X$ hitting $i$ before drifting
away to $\infty$ is strictly smaller than $1$ so that the
$\PP_{i,x}$-expectation over this set in (\ref{31}) equals
$\infty$ (since the integrand tends to $1$ as $t$ tends to
$\infty$) showing that the stopping time (\ref{33}) cannot be
optimal. The proof above shows that the optimality of (\ref{35}) in
this case is obtained through $\tau_n = \tau_{i_n} \wedge
\tau_{f_n}$which play the role of approximate stopping times
(obtained by passing to the limit when $n$ tends to $\infty$ in
(\ref{321}) above). This completes the proof of the theorem.
\end{pf}

\section{The golden ratio rule}\label{sec4}

In this section we show that the minimal solution to (\ref{34})
admits a simple closed-form expression when $X$ is a transient
Bessel process (Theorem~\ref{theor4}). In the case when $X$ is the radial part
of three-dimensional Brownian motion this leads to the golden ratio
rule (Corollary~\ref{coro5}). We also show that $X$ stopped according to the
golden ratio rule has what we refer to as the golden ratio
distribution (Corollary~\ref{coro8}).

In the setting of (\ref{26})--(\ref{28}) consider the optimal
prediction problem (\ref{214}). Recall that this problem is
equivalent to the optimal stopping problem (\ref{221}) which
further extends as (\ref{31}). The main result of this section can
now be stated as follows.
\begin{theorem}\label{theor4}
If $X$ is the $d$-dimensional Bessel
process solving (\ref{26}) with \mbox{$d>2$}, then the optimal stopping
time in (\ref{214}) is given by
%
%
\begin{equation} \label{41}
\tau_* = \inf\{ t \ge0 \vert X_t \ge\lambda I_t \},
\end{equation}
where $\lambda$ is the unique solution to
%
%
\begin{eqnarray}
\label{42}
\lambda^d - (1 + d) \lambda^2 + \frac{4}{4-d} \lambda^{4-d}
- \frac{(d-2)^2}{4-d} &=& 0 \qquad\mbox{if } d \ne4, \\
\label{43}
\lambda^4 - 5 \lambda^2 + 4 \log\lambda+ 4 &=& 0\qquad
\mbox{if } d = 4
\end{eqnarray}
belonging to $(2^{1/(d-2)},\infty)$. The value function (\ref{31})
is given explicitly by
%
%
\begin{eqnarray} \label{44}
V(i,x) &=& \frac{2}{d-2} \biggl[ x^2 \biggl(\frac{1}{2} + \biggl(
\frac{i}{x} \biggr)^{d-2} \biggr) \biggl( \biggl(\frac{\lambda i}{x}
\biggr)^2 - 1 \biggr) \nonumber\\
&&\hspace*{29.7pt}{} - \frac{x^2}{d} \biggl(
\biggl( \frac{\lambda i}{x}\biggr)^d - 1 \biggr) - \frac{2 \lambda^{
4-d}}{d-4} i^2 \biggl( \biggl(\frac{\lambda i}{x}\biggr)^{d-4} - 1
\biggr) \biggr] \qquad\mbox{if } d \ne4 \nonumber\\[-8pt]\\[-8pt]
&=&
\biggl[ x^2 \biggl( \frac{1}{2} + \biggl(\frac{i}{x} \biggr)^2 \biggr) \biggl(
\biggl(\frac{\lambda i}{x}\biggr)^2 - 1 \biggr) \nonumber\\
&&\hspace*{2.6pt}{} - \frac{x^2}{4} \biggl( \biggl( \frac{\lambda i}{x}
\biggr)^4 - 1 \biggr) - 2 i^2 \log\biggl( \frac{
\lambda i}{x} \biggr) \biggr] \qquad\mbox{if } d = 4\nonumber
\end{eqnarray}
for $i \le x \le\lambda i$ and $V(i,x)=0$ for $x \ge\lambda i$
with $i>0$.
\end{theorem}
\begin{pf}
By the result of Theorem~\ref{theor3} we know that the optimal
stopping time $\tau_*$ is given by (\ref{33}) above where the
optimal boundary $f_*$ can be characterized as the minimal solution
to (\ref{34}) staying strictly above the curve $h(i) =
L^{-1}(L(i)/2)$ for $i>0$. Using (\ref{27}) and (\ref{32}) it can
be verified that (\ref{34}) reads as follows:
%
%
\begin{eqnarray} \label{45}
f'(i) &=& \biggl({\frac{d-2}{4-d} \biggl(\frac{f(i)}{i}\biggr) \biggl[ (4
- d) \biggl(\frac{f(i)}{i}\biggr)^{d-2} + (d - 2) \biggl(\frac{f(i)}
{i}\biggr)^{d-4} - 2 \biggr]}\biggr)\nonumber\\
&&\hspace*{0pt}{}\times\biggl({\biggl( \biggl(\frac{f(i)}{i}\biggr)^{d-2} -
1 \biggr) \biggl( \biggl( \frac{f(i)}{i}\biggr)^{d-2} - 2 \biggr)}\biggr)^{-1}
\qquad\mbox{if } d \ne4 \\
&=& \frac{2 ({f(i)}
/{i}) [ ({f(i)}/{i})^2 -2 \log({f(i)}/
{i}) - 1 ]}{( ({f(i)}/{i})^2 - 1 )
( ( {f(i)}/{i})^2 - 2 )} \qquad\mbox{if }
d = 4\nonumber
\end{eqnarray}
and $h(i) = 2^{1/(d-2)} i$ for $i>0$. Hence it is enough to show
that $f_*(i) = \lambda i$ is the minimal solution to (\ref{45})
staying strictly above the curve $h(i) = 2^{1/(d-2)} i$ for $i>0$.

To show that $f_*$ is a solution to (\ref{45}) staying strictly
above $h$, insert $f(i) = \lambda i$ into (\ref{45}) with
$\lambda>0$ to be determined. Multiplying both sides of the
resulting identity by $\lambda^{4-d}$ (to be able to derive the
factorization (\ref{47}) below) it is easy to see that this yields
the equation $F(\lambda)=0$ where we set
%
%
\begin{eqnarray} \label{46}
F(\lambda) &=& \lambda^d - (1 + d) \lambda^2 + \frac{4}{4-d}
\lambda^{4-d} - \frac{(d-2)^2}{4-d} \qquad\mbox{if } d
\ne4 \nonumber\\[-8pt]\\[-8pt]
&=& \lambda^4 - 5 \lambda^2 + 4 \log
\lambda+ 4 \qquad\mbox{if } d = 4\nonumber
\end{eqnarray}
for $\lambda> 0$. After some algebraic manipulations we find that
%
%
\begin{equation} \label{47}
F'(\lambda) = d \lambda^{3-d} \biggl( \lambda^{d-2} - \frac{2}{d}
\biggr) ( \lambda^{d-2} - 2 )
\end{equation}
for $\lambda> 0$ and $d>2$. Hence we see that the equation
$F'(\lambda) = 0$ has two roots $\lambda_0 = (2/d)^{1/(d-2)}$ and
$\lambda_1 = 2^{1/(d-2)}$ where $0 < \lambda_0 < 1 < \lambda_1 <
\infty$. It is easy to check that $F''(\lambda_0) < 0$ and
$F''(\lambda_1) > 0$ showing that $F$ has a local maximum at
$\lambda_0$ and $F$ has a local minimum at $\lambda_1$. Noting that
$F(0+)<0$, $F(1)=0$ and $F(\infty-)=\infty$ this shows that (i) $F$
is strictly increasing on $(0,\lambda_0)$ with $F(0+)<0$ and
\mbox{$F(\lambda_0)>0$}; (ii) $F$ is strictly decreasing on
$(\lambda_0,\lambda_1)$ with $F(1)=0$ and \mbox{$F(\lambda_1)<0$}; and (iii) $F$
is strictly increasing on $(\lambda_1,\infty)$ with
$F(\infty-)=\infty$. It follows therefore that the equation
$F(\lambda)=0$ has exactly three roots $\lambda_0^* < 1 <
\lambda_1^*$ where $\lambda_0^* \in(0,\lambda_0)$ and $\lambda_1^*
\in(\lambda_1,\infty)$. Setting $\lambda= \lambda_1^*$ this shows
that $f_*(i) = \lambda i$ is a solution to (\ref{45}) staying
strictly above the curve $h(i) = 2^{1/(d-2)} i$ for $i>0$ as
claimed.

To show that $f_*$ is the minimal solution satisfying this property,
set $\kappa(i) = f(i)/i$ and note that (\ref{45}) can then be
rewritten as follows:
%
%
\begin{equation} \label{48}
i \kappa'(i) = - \frac{F( \kappa(i) )}{\kappa^{3-d}(i)
( \kappa^{d-2}(i) - 1 ) ( \kappa^{d-2}(i) - 2 )}
\end{equation}
for $i>0$. Since $F(\kappa(i))<0$ for $\kappa(i) \in
(2^{1/(d-2)},\lambda)$ we see from (\ref{48}) that $i \mapsto
\kappa(i)$ is increasing for $\kappa(i) \in(2^{1/(d-2)},\lambda)$.
Noting that (\ref{48}) implies that
%
%
\begin{equation} \label{49}
- \int_{\kappa(i)}^{\kappa(i_0)} \frac{\kappa^{3-d} (\kappa^{d-2}
- 1 ) ( \kappa^{d-2} - 2 )}{F(\kappa)} \,d\kappa=
\int_i^{i_0} \frac{di}{i} = \log\biggl( \frac{i_0}{i} \biggr),
\end{equation}
it follows therefore that the integrand on the left-hand side is
bounded by a constant (not dependent on $i$) as long as
$\kappa(i) \in(2^{1/(d-2)},\lambda)$ for $i \in(0,i_0)$ with any
$i_0>0$ given and fixed. Letting then $i \downarrow0$ in
(\ref{49}) we see that the left-hand side remains bounded while the
right-hand side tends to $\infty$ leading to a contradiction. Noting
that $\kappa(i) \in(2^{1/(d-2)},\lambda)$ if and only if $f(i) \in
(h(i),f_*(i))$, we can therefore conclude that there is no solution
$f$ to (\ref{45}) satisfying $f(i) \in(h(i),f_*(i))$ for $i>0$.
Thus $f_*$ is the minimal solution to (\ref{45}) staying strictly
above $h$ and the proof is complete.
\end{pf}
\begin{corollary}[(The golden ratio rule)]\label{coro5}
If $X$ is the
radial part of three-dimensional Brownian motion, then the optimal
stopping time in (\ref{214}) is given by
%
%
\begin{equation} \label{410}
\tau_* = \inf\biggl\{ t \ge0 \Big\vert\frac{X_t -
I_t}{I_t} \ge\varphi\biggr\},
\end{equation}
where $\varphi= (1 + \sqrt{5})/2 = 1.61 \ldots$ is the golden
%
%
\begin{figure}

\includegraphics{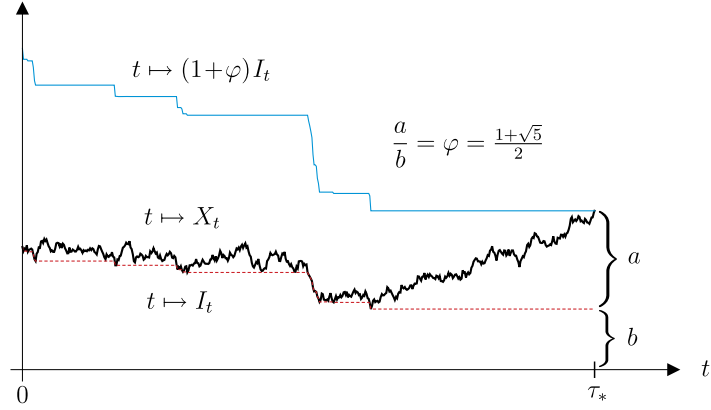}

\caption{The golden ratio rule for the radial part $X$ of
three-dimensional Brownian motion.}\label{fig2}
\end{figure}
ratio (see Figure~\ref{fig2}).
\end{corollary}
\begin{pf}
In this case $d=3$ and equation (\ref{42})
reads
%
%
\begin{equation} \label{411}
\lambda^3 - 4 \lambda^2 + 4 \lambda- 1 = (\lambda- 1) (\lambda^2
- 3\lambda+ 1) =0
\end{equation}
for $\lambda>0$. Solving the latter quadratic equation and choosing
the root strictly greater than $1$, we find that $\lambda= 1 +
\varphi$ where $\varphi= (1 + \sqrt{5})/2 = 1.61 \ldots$ is the
golden ratio. The optimality of (\ref{410}) then follows from
(\ref{41}) and the proof is complete.
\end{pf}

Returning to the result of Theorem~\ref{theor3} above we now determine the law
of the transient diffusion process $X$ stopped at the optimal
stopping time $\tau_*$ (for related results on the Skorokhod
embedding problem see~\cite{RY}, pages~269--277, and the references
therein).
\begin{proposition}\label{prop6}
In the setting of Theorem~\ref{theor3} we have
%
%
\begin{equation} \label{412}
\PP_{ x}( X_{\tau_*} \le y ) = \exp\biggl( -
\int_{f_*^{-1}(y)}^x \frac{dL(z)}{L(f_*(z)) - L(z)} \biggr)
\end{equation}
for $0 < y \le f_*(x)$ with $x>0$.
\end{proposition}
\begin{pf}
Note that
%
%
\begin{eqnarray} \label{413}
\tau_* &=& \inf\{ t \ge0 \vert X_t \ge f_*(I_t) \} \nonumber\\
&=& \inf\{ t \ge0 \vert L(X_t) \ge(L \circ f_* \circ
L^{-1}) (L(I_t)) \} \\
&=& \inf\{ t \ge0 \vert
X_t^L \ge f_*^L(I_t^L) \},\nonumber
\end{eqnarray}
where we set $X_t^L = L(X_t)$, $f_*^L = L \circ f_* \circ L^{-1}$
and $I_t^L = L(I_t) = \inf_{ 0 \le s \le t} L(X_s) = \inf_{ 0 \le
s \le t} X_s^L$ for $t \ge0$. Let $x>0$ be given and fixed. For $j
\le L(x)$ set $G(j) = \int_{-\infty}^j g(k) \,dk$ where $g\dvtx
(-\infty,L(x)] \rightarrow\R$ is a continuously differentiable
function with bounded support. Using the fact that $dI_t^L=0$ when
$X_t^L \ne I_t^L$ it is easily verified by It\^o's formula that the
process $M^L = (M^L_t)_{t \ge0}$ defined by
%
%
\begin{equation} \label{414}
M_t^L = G(I_t^L) + (X_t^L - I_t^L) G'(I_t^L)
\end{equation}
is a continuous local martingale. Moreover, since $G'=g$ is
continuous and has bounded support, we see that $M^L$ is bounded and
therefore uniformly integrable. By the optional sampling theorem we
thus find that
%
%
\begin{eqnarray} \label{415}
G^L(x) &=& \EE_x M_0^L = \EE_x M_{\tau_*}^L = \EE_x G(I_{\tau_*}^L) +
\EE_x [ (X_{\tau_*}^L - I_{\tau_*}^L) G'(I_{\tau_*}^L) ]
\nonumber\\
&=& \int_{-\infty}^{L(x)} G(j) \,dF(j) + \EE_x \bigl[
\bigl(f_*^L (I_{\tau_*}^L) - I_{\tau_*}^L\bigr) G'(I_{\tau_*}^L)
\bigr]\nonumber\\
&=& G(j) F(j) \vert_{-\infty}^{L(x)} - \int_{-\infty
}^{L(x)} F(j) \,dG(j)\nonumber\\[-8pt]\\[-8pt]
&&{} + \int_{-\infty}^{L(x)} \bigl(f_*^L(j) - j\bigr) G'(j)
\,dF(j) \nonumber\\
&=& G^L(x) - \int_{-\infty}^{L(x)} F(j) g(j)
\,dj\nonumber\\
&&{} + \int_{-\infty}^{L(x)} \bigl(f_*^L(j) - j\bigr) g(j) \,dF(j),\nonumber
\end{eqnarray}
where we set $G^L(x) = G(L(x))$ and $F$ denotes the distribution
function of $I_{\tau_*}^L$ under $\PP_{ x}$. Since (\ref{415})
holds for all functions $g$ of this kind, it follows that
%
%
\begin{equation} \label{416}
F'(j) = \frac{F(j)}{f_*^L(j) - j}
\end{equation}
for $j<L(x)$ with $F(L(x))=1$. Solving (\ref{416}) under this
boundary condition we find that
%
%
\begin{equation} \label{417}
F(j) = \exp\biggl( - \int_j^{L(x)} \frac{dk}{f_*^L(k) - k}
\biggr)
\end{equation}
for $j \le L(x)$. Recalling that $f_*^L = L \circ f_* \circ L^{-1}$
and substituting $k=L(z)$ it follows that
%
%
\begin{eqnarray} \label{418}\qquad
\PP_{ x}(I_{\tau_*} \le i) &=& \PP_{ x}\bigl(L(I_{\tau_*}) \le L(i)\bigr) =
\PP_{ x}\bigl(I_{\tau_*}^L \le L(i)\bigr)\nonumber\\
&=& F(L(i))
= \exp
\biggl( - \int_{L(i)}^{L(x)} \frac{dk}{f_*^L(k) - k} \biggr)\\
&=& \exp\biggl( - \int_i^x \frac{dL(z)}{L(f_*(z)) - L(z)} \biggr)\nonumber
\end{eqnarray}
for $i \le x$ in $(0,\infty)$. Hence we find that
%
%
\begin{eqnarray} \label{419}
\PP_{ x}(X_{\tau_*} \le y) &=& \PP_{ x}\bigl(f_*(I_{\tau_*}) \le y\bigr)
= \PP_{ x}\bigl(I_{\tau_*} \le f_*^{-1}(y)\bigr) \nonumber\\[-8pt]\\[-8pt]
&=& \exp\biggl( -
\int_{f_*^{-1}(y)}^x \frac{dL(z)}{L(f_*(z)) - L(z)} \biggr)\nonumber
\end{eqnarray}
for $0 < y \le f_*(x)$ with $x>0$. This completes the proof.\vadjust{\goodbreak}
\end{pf}

Specializing this result to the $d$-dimensional Bessel process
$X$ of Theorem~\ref{theor4} we obtain the following consequence.
\begin{corollary}\label{coro7}
In the setting of Theorem~\ref{theor4} we have
%
%
\begin{equation} \label{420}
\PP_{ x}( X_{\tau_*} \le y ) = \biggl( \frac{y}{\lambda
x} \biggr)^{({d-2})/({1-(1/\lambda)^{d-2}})}
\end{equation}
for $0 < y \le\lambda x$ with $x>0$.
\end{corollary}
\begin{pf}
In this case $f_*(i) = \lambda i$ for $i>0$ where
$\lambda$ is the unique solution to either (\ref{42}) when $d \ne
4$ or (\ref{43}) when $d=4$ and $L$ is given by (\ref{27}).
Inserting these expressions into the right-hand side of (\ref{412})
it is easily verified that this yields (\ref{420}).
\end{pf}

Specializing this further to the radial part $X$ of
three-dimensional Brownian motion in Corollary~\ref{coro5} we obtain the
following conclusion.
\begin{corollary}[(The golden ratio distribution)]\label{coro8}
In the
setting of Corollary~\ref{coro5} we have
%
%
\begin{equation} \label{421}
\PP_{ x}( X_{\tau_*} \le y ) = \biggl( \frac{y}{(1 + \varphi)
x} \biggr)^\varphi
\end{equation}
for $0 < y \le(1 + \varphi) x$ with $x>0$.
\end{corollary}
\begin{pf}
In this case $d=3$ and $\lambda= 1+\varphi$ so that
$(d - 2)/(1 - (1/\lambda)^{d-2}) = 1/(1 - 1/(1 + \varphi)) = (1 +
\varphi)/\varphi= \varphi^2/\varphi= \varphi$. Hence we see that
(\ref{420}) reduces to (\ref{421}).
\end{pf}

Note from (\ref{421}) that the density function of $X_{\tau_*}$
under $\PP_{ x}$ is given by
%
%
\begin{equation} \label{422}
f_{X_{\tau_*}}(y) = \frac{\varphi}{((1 + \varphi) x)^\varphi}
y^{\varphi-1}
\end{equation}
for $0 < y < (1 + \varphi) x$ with $x>0$ and equals zero
otherwise. We refer to (\ref{421}) and (\ref{422}) as the \textit{golden
ratio distribution}. It is easy to see that
%
%
\begin{equation} \label{423}
\EE_x X_{\tau_*} = \varphi x
\end{equation}
for $x>0$. The fact that this number is strictly greater than $x$
(the initial point corresponding to stopping at once) is not
surprising since $X$ is a submartingale. It needs to be recalled
moreover that the aim of applying the golden ratio rule $\tau_*$ is
to be as close as possible to the time $\theta$ at which the
ultimate minimum is attained. We will see in the next section that
the golden ratio distribution provides insight as to what extent the
golden ratio rule has the power of capturing the ultimate maximum of
a strict local martingale.\looseness=1

\section{Applications in optimal trading}\label{sec5}

In this section we present some applications of the previous results
in problems of optimal trading. We also outline some remarkable
connections between such problems and the practice of technical
analysis. These applications and connections rest on three basic
ingredients that we describe first.\vspace*{8pt}

1. \textit{Fibonnaci retracement.} We begin by explaining a few
technical terms from the field of applied finance. \textit{Technical
analysis} is a financial term used to describe methods and
techniques for forecasting the direction of asset prices through the
study of past market data (primarily prices themselves plus the
volume of their trade). \textit{Support} and \textit{resistance} are
concepts in technical analysis associated with the expectation that
the movement of the asset price will tend to cease and reverse its
trend of decrease/increase at certain predetermined price levels. A
\textit{support}/\textit{resistance level} is a price level at which the
price will tend to find support/resistance when moving down/up. This
means that the price is more likely to bounce off this level rather
than break through it. One may also think of these levels as turning
points of the prices. \textit{Fibonacci retracement} is a method of
technical analysis for determining support and resistance levels.
The name comes after its use of \textit{Fibonacci numbers} $F_{n+1} =
F_n + F_{n-1}$ for $n \ge1$ with $F_0=0$ and $F_1=1$. Fibonacci
retracement is based on the idea that after reversing the trend at a
support/resistance level, the price will retrace a predictable
portion of the past downward/upward move by advancing in the
opposite direction until finding a new resistance/support level,
after which it will return to the initial trend of moving
downwards/upwards. Fibonacci retracement is created by taking two
extreme points on a chart showing the asset price as a function of
time and dividing the vertical distance between them by the key
Fibonacci ratios ranging from $0\%$ (start of the retracement) to
$100\%$ (end of the retracement representing a complete reversal to
the original trend). The other key Fibonacci ratios are $23.6\%$
(shallow retracement), $38.2\%$ (moderate retracement) and $61.8\%$
(golden retracement). They are obtained by formulas $(F_n/F_{n+3})
\times100 \approx\varphi^{-3} \times100$, $(F_n/F_{n+2}) \times
100 \approx\varphi^{-2} \times100$ and $(F_n/F_{n+1}) \times100
\approx\varphi^{-1} \times100$, respectively (see the next
paragraph). These retracement levels serve as alert points for a
potential reversal at which traders may employ other methods of
technical analysis to identify and confirm a reversal. Despite its
widespread use in technical analysis of asset prices, there appears
to be no (rigorous) explanation of any kind as to why the Fibonacci
ratios should be used to this effect. We will show below that the
golden ratio rule derived in the previous section offers a rigorous
optimality argument for the choice of the golden retracement
($61.8\%$). To our knowledge this is the first time that such an
argument has been found/given in the literature.\vadjust{\goodbreak}

2. \textit{Golden ratio and Fibonacci numbers.} The link between the
two is well known and is expressed by Binet's formula
%
%
\begin{equation} \label{51}
F_n = \frac{\varphi^n - \psi^n}{\varphi- \psi} = \frac{\varphi^n
- \psi^n}{\sqrt{5}},
\end{equation}
where $\varphi= (1 + \sqrt{5})/2$ and $\psi= (1 - \sqrt{5})/2 =
1 - \varphi= -1/\varphi$. It follows that
%
%
\begin{equation} \label{52}
\lim_{n \rightarrow\infty} \frac{F_{n+1}}{F_n} = \varphi.
\end{equation}
This fact is used in the description of Fibonacci retracement above.\vspace*{8pt}

3. \textit{The CEV model.} One of the simplest/tractable models for
asset price movements that is capable of reproducing the implied
volatility smile/frown effect and the (inverse) leverage effect
(both observed in the empirical data) is the \textit{Constant
Elasticity of Variance} (CEV) model in which the (nonnegative)
asset price process $Z=(Z_t)_{t \ge0}$ solves
%
%
\begin{equation} \label{53}
dZ_t = \mu Z_t \,dt + \sigma Z_t^{1+\beta} \,dB_t,
\end{equation}
where $\mu\in\R$ is the appreciation rate, $\sigma>0$ is the
volatility coefficient, and $\beta\in\R$ is the elasticity
parameter. If $\beta=0$ then $Z$ is a geometric Brownian motion
which was initially considered in~\cite{Os} and~\cite{Sa}. For
$\beta\ne0$ this model was firstly considered in~\cite{Cox} for
$\beta<0$ and then in~\cite{EM} for $\beta>0$. Due to its predictive
power and tractability, the CEV model is widely used by
practitioners in the financial industry, especially for modeling
prices of equities and commodities. If $\beta<0$ then the model
embodies the \textit{leverage effect} (commonly observed in equity
markets) where the volatility of the asset price increases as its
price decreases. If $\beta>0$ then the model embodies the
\textit{inverse leverage effect} (often observed in commodity markets)
where the volatility of the asset price increases when its price
increases. For example, it is reported in~\cite{GS} that the
elasticity coefficient $\beta$ for Gold on the London Bullion Market
in the period from 2000 to 2007 was approximately 0.49. Similar
elasticity coefficients have also been observed for other precious
metals (such as Copper for instance).

In the remainder of this section we focus on the case when $\mu=0$
and \mbox{$\beta>0$}. It is well known (cf.~\cite{EM}) that $Z$ solving
(\ref{53}) is a strict local martingale (a local martingale which
is not a true martingale) in this case due to the fact that $t
\mapsto\EE_z(Z_t)$ is strictly decreasing on $\R_+$ for any $z>0$.
This also implies that $Z$ does not admit an equivalent martingale
measure so that the CEV model may admit arbitrage opportunities. One
way of looking at the models of this type is to associate them with
asset price bubbles (see~\cite{HLW}). After soaring to a finite
ultimate maximum (bubble) at a finite time, the asset price will
tend to zero as time goes to infinity, and the central question for
a holder of the asset becomes when to sell so as to be as close as
possible to the time at which the ultimate maximum is attained. More
precisely, introducing the running maximum process $S=(S_t)_{t \ge
0}$ associated with $Z$ by setting
%
%
\begin{equation} \label{54}
S_t = \sup_{0 \le s \le t} Z_s
\end{equation}
and recalling that $Z_t \rightarrow0$ as $t \rightarrow\infty$, we
see that the ultimate supremum $S_\infty= \sup_{ t \ge0} Z_t$ is
attained at some random time $\theta$ in the sense that
%
%
\begin{equation} \label{55}
Z_\theta= S_\infty
\end{equation}
with $\PP_{ z}$-probability one for $z>0$ given and fixed. The
optimal selling problem addressed above then becomes the optimal
prediction problem
%
%
\begin{equation} \label{56}
V(z) = \inf_\tau\EE_z ( \vert\theta- \tau\vert-
\theta),
\end{equation}
where the infimum is taken over all stopping times $\tau$ of $Z$
(with finite mean) and $z>0$ is given and fixed. We will now show
that due to the well-known connection between CEV and Bessel
processes (dating back to similar transformations in~\cite{Cox} and
\cite{EM}) the problem (\ref{56}) can be reduced to the problem
(\ref{214}) solved above.\vspace*{8pt}

\begin{figure}

\includegraphics{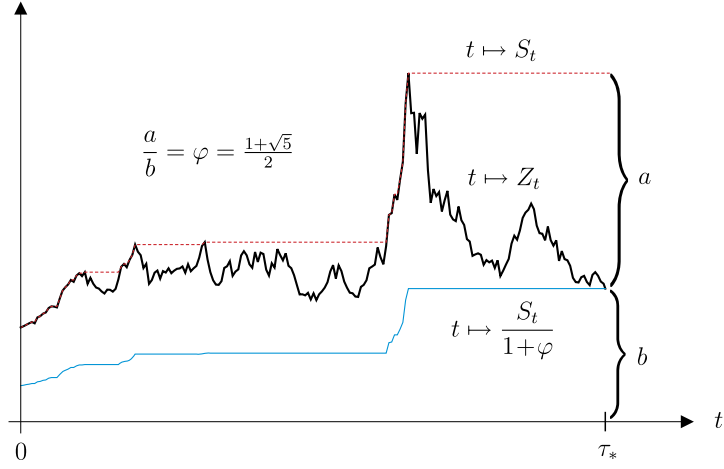}

\caption{The golden ratio rule for the CEV process $Z =
1/X$ where $X$ is the radial part of three-dimensional Brownian
motion. Note the presence of a bubble and its relation to the golden
ratio.}\label{fig3}
\end{figure}

4. \textit{The golden ratio rule for the CEV process.} For $d>2$
given and fixed consider the $d$-dimensional Bessel process $X$
solving (\ref{26}) under $\PP_{ x}$ with $x>0$. Recall that the
scale function $L$ is given by (\ref{27}) and set $K(x) =
-c_\sigma L(x)$ for $x>0$ with $c_\sigma>0$ given and fixed.
Then the process $Z=K(X)$ defined by
%
%
\begin{equation} \label{57}
Z_t = K(X_t) = \frac{c_\sigma}{X_t^{d-2}}
\end{equation}
is on natural scale and It\^o's formula shows that $Z$ solves
%
%
\begin{equation} \label{58}
dZ_t = \sigma Z_t^{1+{1}/({d-2})} \,d\tilde B_t,
\end{equation}
where $\sigma= (d - 2)/c_\sigma^{1/(d-2)}$ and $\tilde B = -B$ is
a standard Brownian motion. Note that equation (\ref{58})
coincides with equation (\ref{53}) for $\mu=0$ and $\beta=
1/(d - 2)$. By the uniqueness in law for this equation (among
positive solutions) it follows that $Z=K(X)$ is a CEV process. From
the properties of $X$ it follows that after starting at $z=K(x)>0$,
the process $Z$ stays strictly positive (without exploding at a
finite time) and $Z_t \rightarrow0$ with $\PP_{ z}$-probability
one as $t \rightarrow\infty$. This shows that $\theta$ in
(\ref{55}) is well defined. Moreover, due to the reciprocal
relationship (\ref{57}) we see that the time of the ultimate
maximum $\theta$ for $Z$ in (\ref{55}) coincides with the time of
the ultimate minimum $\theta$ for $X$ in (\ref{210}), and hence the
problem (\ref{56}) has the same solution as the problem
(\ref{214}) (note also that the natural filtrations of $Z$ and $X$
coincide so that $\tau$ is a stopping time of $Z$ if and only if
$\tau$ is a stopping time of $X$). Since $X_t \ge\lambda I_t$ if
and only if $Z_t/c_\sigma= X_t^{2-d} \le\lambda^{2-d}
I_t^{2-d} = \lambda^{2-d} S_t/c_\sigma$ it follows from
(\ref{41}) in Theorem~\ref{theor4} that the optimal stopping time in
(\ref{56}) is given by
%
%
\begin{equation} \label{59}
\tau_* = \inf\{ t \ge0 \vert S_t \ge\lambda^{d-2} Z_t \},
\end{equation}
where $\lambda$ is the unique solution to either (\ref{42}) or
(\ref{43}) belonging to $(2^{1/(d-2)},\infty)$. In particular, if
$d=3$ then we know from (\ref{411}) that $\lambda= 1 + \varphi$
so that (\ref{59}) reads
%
%
\begin{equation} \label{510}
\tau_* = \inf\biggl\{ t \ge0 \Big\vert\frac{S_t -
Z_t}{Z_t} \ge\varphi\biggr\}.
\end{equation}
This is the \textit{golden ratio rule} for the CEV process $Z = 1/X$
where $X$ is the radial part of three-dimensional Brownian motion (see
Figure~\ref{fig3}).


%
%


To relate the golden ratio rule (\ref{510}) to Fibonacci
retracement discussed above, let $a = S_{\tau_*} - Z_{\tau_*}$
denote the larger quantity and let $b=Z_{\tau_*}$ denote the smaller
quantity in the golden ratio rule. To determine the percentage of
$a$ in $a + b$ we need to calculate
%
%
\begin{eqnarray} \label{511}
\frac{a}{a + b} &=& \frac{S_{\tau_*} - Z_{\tau_*}}{S_{\tau_*}} = 1
- \frac{Z_{\tau_*}}{S_{\tau_*}} = 1 - \frac{1}{1 + \varphi}\nonumber\\[-8pt]\\[-8pt]
&=&
\frac{\varphi}{1 + \varphi} = \frac{1}{\varphi}.\nonumber
\end{eqnarray}
Multiplying this expression by $100$ gives $61.8\%$ and this is
exactly the \textit{golden retracement} discussed above. In view of
the optimality of (\ref{510}) in (\ref{56}) we see that the golden
retracement of $61.6\%$ for the CEV process $Z = 1/X$ (starting
close to zero) where $X$ is the radial part of three-dimensional
Brownian motion can be seen as a rational support level (in the
sense that rational investors who aim at selling the asset at the
time of the ultimate maximum will sell the asset at the time of the
golden retracement, and therefore the asset price could be expected
to raise afterwards). To our knowledge this is the first time that
such a rational optimality argument for the golden retracement has
been established.



\printaddresses

\end{document}